\pdfoutput=1
\documentclass[11pt,a4paper,twoside]{article}

\usepackage[T1]{fontenc}
\usepackage[utf8]{inputenc}
\usepackage{latexsym}
\usepackage{amsfonts,amsmath,amssymb,amsthm}
\usepackage{fullpage}
\usepackage{xcolor,bm,bbm}
\usepackage{array}
\usepackage{mathrsfs}
\usepackage{mathtools}
\usepackage{fourier}
\usepackage{graphicx}
\usepackage{hyperref}

\renewcommand{\P}{\mathbb P}

\newcommand{\E}{\mathbb E}
\newcommand{\V}{\mathbb V}

\newcommand{\R}{\mathbb R}
\newcommand{\N}{\mathbb N}

\renewcommand{\O}{\mathcal O}
\newcommand{\dint}{\,\mathrm{d}}
\DeclareMathOperator{\Var}{Var}
\DeclareMathOperator{\Unif}{Unif}
\DeclareMathOperator{\vol}{vol}

\DeclareMathOperator{\Tr}{Tr}

\newtheorem{thm}{Theorem}[section]

\newtheorem{cor}[thm]{Corollary}
\newtheorem{lemma}[thm]{Lemma}
\newtheorem{df}[thm]{Definition}
\newtheorem{proposition}[thm]{Proposition}

\newtheorem{thmalpha}{Theorem}

\theoremstyle{definition}

\newtheorem{rmk}[thm]{Remark}

\allowdisplaybreaks
\hypersetup{
  colorlinks=true,
  linkcolor=black,
  citecolor=black,
  urlcolor=black,
  pdftitle={Large and Moderate Deviations for Entries of Orthogonal Matrices and the Stiefel Manifold},
  pdfauthor={Philipp Tuchel}
}

\begin{document}

\title{\textbf{Large and Moderate Deviations for Entries\\
of Orthogonal Matrices and the Stiefel Manifold}}
\author{Philipp Tuchel\thanks{Faculty of Mathematics, Ruhr University Bochum, Universit\"atsstra\ss e 150, 44780 Bochum, Germany. Email: \href{mailto:philipp.tuchel@rub.de}{\nolinkurl{philipp.tuchel@rub.de}}.}}
\date{}
\maketitle

\begin{abstract}
\small
Let $A_N$ be distributed according to the Haar probability measure on the orthogonal group $\O(N)$ for each $N\in\N$. It is well-known that the upper left $m_N\times k_N$ block of $\sqrt{N}A_N$ with $m_Nk_N = o(N)$ converges in total variation distance to a matrix of same size consisting of i.i.d. standard normal entries as $N\to\infty$. In this work, we characterize this convergence on the scale of large deviations. More precisely, we show that under the same condition $m_Nk_N = o(N)$ the empirical measure of entries of this block satisfies a large deviation principle with speed $m_Nk_N$ and rate function given by the relative entropy with respect to the standard normal distribution. Further, we complement the large deviation principle (LDP) obtained by Kabluchko and Prochno in [Large deviations for random matrices in the orthogonal group and Stiefel manifold with applications to random projections of product distributions, Annales de l'Institut Henri Poincar\'e. 60 (2024), 990 -- 1024] for the entire matrix $A_N$ with a moderate deviation principle (MDP). Concretely, we show an MDP for the sequence of matrices $\beta_N A_N$ in the product topology, where $\beta_N\to\infty$ is a sequence of real numbers such that $\beta_N = o(\sqrt{N})$. Here, in contrast to the LDP, the Gaussian behavior of the entries is reflected in the rate function.

\medspace
\end{abstract}

\section{Introduction}
\label{C:sec:introduction}

In this work, we consider the following sets of matrices that carry invariant probability measures. For integers \(m \leq N\), the \emph{Stiefel manifold} \(\mathbb{V}_{m,N}\) is the set of all \(m \times N\) real matrices with orthonormal rows:
	\begin{align}
		\mathbb{V}_{m,N} \; := \; \Bigl\{\, V \in \mathbb{R}^{m\times N} \;:\; V\,V^\top \;=\; \operatorname{Id}_m \Bigr\}.
	\end{align}
	Here, $\operatorname{Id}_m$ denotes the $m$-dimensional identity matrix. When \(m = N\), this coincides with the \emph{orthogonal group} \(\O(N)\), which is the set of all \(N \times N\) real matrices satisfying $A A^\top = A^\top A = \operatorname{Id}_N$. For these sets (which have a group structure only in special cases, e.g. $m=N$), there is a unique invariant probability measure, called Haar measure, that we denote by $\Unif(\V_{m,N})$ and $\Unif(\O(N))$ in case of $m= N$. Its characterizing property is that for any fixed $A\in \O(m)$ and $B\in \O(N)$, the distributions of $A V B$ and $V$ coincide for $V$ being distributed as $\Unif(\V_{m,N})$, denoted by $V\sim \Unif(\V_{m,N})$. In \cite{james1954normal} James connected these measures to known matrix distributions such as the Wishart distribution. For a modern, comprehensive treatment of random matrix theory on classical compact groups, we refer the reader to the monograph by E.\ Meckes \cite{meckes2019random}. The Haar measure provides a natural framework for the study of the typical behavior of such sets of matrices. In particular, the distribution of the entries of random orthogonal matrices in high dimensions has attracted considerable attention and is an active research area in probability and statistics today. Roughly speaking, as each row (and column) of a matrix in $\O(N)$ has norm $1$, a typical entry of a random matrix in $\O(N)$ is of order $\frac{1}{\sqrt{N}}$. Over a century ago Poincaré and later Maxwell and Borel (see \cite{borel1914introduction} and \cite[Chapter 6]{diaconis1987dozen} for historical context) observed that any individual coordinate of a sphere, rescaled by $\sqrt{N}$, is approximately Gaussian for large $N$. This classical result is a precursor to the study of collective convergence of coordinates to i.i.d. Gaussian random variables. In \cite{diaconis1987dozen}, Diaconis and Freedman proved that $p_N$ entries of a random uniform vector on the sphere in $\R^N$ converge jointly in total variation to $p_N$ i.i.d. Gaussian random variables provided $p_N = o(N)$. A notable result in this direction is by Jiang who proved a threshold behavior for the more general case of the orthogonal group $\O(N)$ (see \cite{jiang2006many}). More precisely, let $A_N\sim \Unif(\O(N))$ and let $Z_N = A_N^{m_N, k_N}$ be the $m_N \times k_N$ upper-left block of $A_N$. Then, the joint distribution of the $m_N k_N$ entries of $\sqrt{N} Z_N$ converges in total variation to the joint distribution of $m_N k_N$ independent standard normal random variables as $N\to\infty$ if $m_N = o(\sqrt{N})$ and $k_N = o(\sqrt{N})$. On the other hand, he proved that if $m_N = \lfloor x\sqrt{N} \rfloor$ and $k_N = \lfloor y\sqrt{N} \rfloor$ for fixed $x, y > 0$, this convergence does not hold. In recent years, more progress has been made in this direction. For instance, in \cite{DistancesJiang} Jiang and Ma proved that the convergence even holds under the condition $m_N k_N = o(N)$ but not if $m_Nk_N = cN$ for some constant $c>0$. Moreover, the convergence to Gaussian entries under the regime $m_Nk_N=o(N)$ holds not only in total variation but also for more distances including the Kullback-Leibler distance and the Hellinger distance. Notably, in the Euclidean distance, the convergence holds if $m_Nk_N^2 = o(N)$ but not if $m_N k_N^2 = cN$ for a positive constant $c$. Recently, in \cite[Theorem 1]{chen2025fisher}, the authors provide the explicit convergence rate $\frac{m_N^2k_N(k_N+1)}{4N^2}(1-o(1))$ for the Fisher information between $\sqrt{N}Z_N$ and an $m_N\times k_N$ block of i.i.d. standard Gaussians.  \\
	A framework for a more quantitative study of the asymptotic behavior beyond the scale of Gaussian approximation is the theory of large deviations. Originated in the work of Cramér on sums of i.i.d. random variables and formalized by Varadhan and others in the 1960s, large deviation principles (LDPs) have become a cornerstone in modern probability theory (see \cite{cramer1938nouveau}, \cite{donsker1975asymptotic}, \cite{varadhan1966asymptotic} and for a book \cite{dembo2009techniques}). In rough terms, LDPs characterize probabilities of rare events on an exponential scale in a quantitative way including the speed of convergence and a so-called rate function that captures the, often subtle, dependence on the event (see Definition \ref{def:LDP_main} for a precise definition). Classical results include Cramér's LDP for sample means of random variables, Sanov's theorem for random empirical measures and Bryc's inverse Varadhan lemma (see \cite[Chapters 4 and 6]{dembo2009techniques}). In cases where, for some sequence of random variables $X_N$, the sequence $X_N$ satisfies an LDP and $\sqrt{N}X_N$ satisfies a central limit theorem (CLT), often the sequence $\beta_NX_N$, for some $\beta_N\to\infty$, $\beta_N = o(\sqrt{N})$, also satisfies an LDP with a certain speed and a rate function that reflects the CLT behavior. Such LDPs are also called moderate deviation principles (MDPs) and provide convergence behavior on this intermediate scale. \\
	In the context of random matrices, large deviation results have been mostly studied in the context of spectral properties, such as the empirical eigenvalue distributions of Wigner and Wishart matrices (see for instance \cite{anderson2010introduction}, \cite{guionnet2004large}, \cite{wigner1958distribution}, \cite{wigner1993characteristic} and \cite{wishart1928generalised}). \\
	In contrast, large deviations for the entries of such matrices are a more recent pursuit. One result in this direction is given by Kim and Ramanan (\cite[Theorem 2.8]{kim2021large}) who showed that the empirical measure of scaled columns of an element from the Stiefel manifold satisfies an LDP. More precisely, $\frac{1}{N} \sum_{i=1}^{N} \delta_{\sqrt{N} v_i}$, where $V_N = (v_1, \ldots, v_N) \sim \Unif(\V_{m,N})$, satisfies an LDP in the $q$-Wasserstein topology for $q\in (0,2)$ (for background, see \cite{panaretos2020invitation}) with speed $N$ and good (strictly) convex rate function $I$ given by 
	\begin{align}\label{eq:ratefunctionstiefelkim}
	I(\nu) = \begin{cases}
		H(\nu| \gamma^{\otimes m}) +\frac{1}{2}\mathrm{Tr}(\operatorname{Id}_m - \mathcal{C}(\nu)) & \text{if } \operatorname{Id}_m-\mathcal{C}(\nu) \text{ is positive semidefinite}, \\
		+\infty & : \text{otherwise,}
	\end{cases}
	\end{align}
	Here $H$ is the relative entropy, $\gamma^{\otimes m}$ is the Gaussian product measure on $\R^m$ and $\mathcal{C}(\nu)$ is the covariance operator of $\nu$. Note that the rate function is, in this case, not just the relative entropy which would be expected from a purely Gaussian limit. This statement is a core element in the proof of their results about large deviations for random projection of certain high-dimensional random vectors (\cite[Theorem 2.4 and Theorem 2.6]{kim2021large}), illustrating the significance of these studies in asymptotic geometric analysis. Random projections of vectors or convex bodies are a common application of random elements from the Stiefel manifold and an active research area (see \cite{Prochno2024} for a recent survey and the references therein). Another related result was obtained in \cite{KLR2022}. The authors consider, for a certain class of random vectors $X_N$ with values in $\R^N$, the empirical measure of the projection coordinates:
	$
	\frac{1}{m_N} \sum_{i=1}^{m_N} \delta_{(A_N X_N)_i}.
	$
	Here, $A_N$ is the upper left $m_N \times N$ block of a random orthogonal matrix and $A_NX_N$ can be seen as the projection of $X_N$ onto a random subspace. They prove an LDP for this sequence with different regimes for $m_N$ (constant, sub-linear and linear). A key ingredient in this proof is an LDP for the empirical measure of $m_N$ coordinates of a uniform random element from the sphere in $\R^N$, i.e. for $\frac{1}{m_N} \sum_{i=1}^{m_N} \delta_{\sqrt{N} X_i}$, where $(X_1,\ldots, X_N)\sim \Unif(\V_{1,N})$. This gives another point of view on the central limit theorem for $m_N$ coordinates of a uniform random element from the sphere in $\R^N$ obtained by Diaconis and Freedman in \cite{diaconis1987dozen} in the context of large deviations and empirical measures. \\
    In addition to the limiting behavior of the entries scaled by $\sqrt{N}$, there has been some progress on the limiting behavior of the unscaled entries of $A_N\sim \Unif(\O(N))$. Kabluchko and Prochno \cite{KabluchkoLarge} (see also Kabluchko, Prochno and Thäle \cite{kabluchko2021new}) quantified the rate of this convergence on the scale of large deviations motivated by applications in random projections of product distributions. They showed (\cite[Theorem B]{KabluchkoLarge}) that the sequence $A_N$ (embedded in \([-1,1]^{\infty\times\infty}\) by padding with zeros) satisfies, for $N\to\infty$, an LDP with speed $N$ and good convex rate function $I:[-1,1]^{\infty\times\infty}\to [0,\infty]$, given by
	\begin{align}\label{ratefunctionorthogonal}
		I(T):= \begin{cases}
			-\frac{1}{2} \log \det \left(\operatorname{Id}_\infty - T T^\top\right) & : T \in \mathscr{S}_2 \text{ and } \|T T^\top\|_{op} < 1, \\
			+\infty & : \text{otherwise.}
		\end{cases}
	\end{align}
	where $\|\cdot\|_{op}$ denotes the operator norm and $\mathscr{S}_2$ is the Hilbert space of Hilbert-Schmidt operators on $\ell_2$: $\mathscr{S}_2 := \{T = (t_{i,j})_{i,j \in \mathbb{N}} \in \mathbb{R}^{\infty \times \infty} : \sum_{i,j \in \mathbb{N}} t_{i,j}^2 < \infty \}$. Here, the determinant is defined as 
  $$
  \det \left(\operatorname{Id}_\infty - T T^\top\right) := \lim_{k \to \infty} \det \left(\operatorname{Id}_k - T_k T_k^\top\right),
  $$
  with $T_k$ being the top-left $k\times k$ block of $T$. Since $T\in\mathscr{S}_2$, the operator $TT^\top$ is trace class, and the determinant is well defined and positive under $\|TT^\top\|_{op}<1$. Similarly, according to \cite[Theorem A]{KabluchkoLarge}, for fixed $m\in\N$ and $V_N \sim \Unif(\V_{m,N})$, the sequence $V_N$ satisfies an LDP with speed $N$ and a good convex rate function $J:[-1,1]^{m\times\infty}\to [0,\infty]$ that has a similar structure to the rate function $I$ above. \\

	The purpose of this work is twofold. First, we characterize the normal approximation of an $m_N\times k_N$ block in a new way, by an LDP for the empirical measure of the scaled entries of a random orthogonal matrix, where the Gaussian behavior is captured in the relative entropy in the rate function. In fact, as a consequence we find that the empirical measure converges, under certain assumptions, almost surely to the standard Gaussian measure in the weak topology. This LDP complements and extends the empirical measure LDPs by Kim, Liao and Ramanan (\cite{kim2021large},\cite{KLR2022}) and gives a quantitative interpretation of the Gaussian approximation on the scale $m_Nk_N = o(N)$. At the heart of the proof lies a new uniform local limit result for densities, which might be of independent interest. Next, we examine the asymptotic behavior of the entire matrix block not on the large but on the moderate deviation scale. This means that for any sequence $\beta_N\to\infty$ and $\beta_N = o(\sqrt{N})$, we show that the sequence $\beta_N A_N$ satisfies an MDP with speed $\frac{N}{\beta_N^2}$ and a rate function that reflects the Gaussian limit of the individual entries, which is not seen in the unscaled regime. This closes the gap between the matrix scaled by $\sqrt{N}$, studied by Jiang and Ma in \cite{DistancesJiang}, and the unscaled matrix, studied by Kabluchko and Prochno in \cite{KabluchkoLarge}. 

\section{Main results}
\label{C:sec:results}
\setcounter{thmalpha}{0}

Our first result is a large deviation principle (LDP) for the empirical measure of the scaled entries of a random orthogonal matrix. The rate function reflects the Gaussian behavior of this empirical measure. 

	\begin{thmalpha}[Empirical Measure LDP for Scaled Entries]\label{thm:empirical_main} 
		Let $X_N = A_N^{m_N, k_N}$ be the $m_N\times k_N$ upper left block of $A_N$ where $A_N$ is Haar-distributed on $\O(N)$. Assume that $m_N, k_N$ are sequences of natural numbers such that $p_N := m_N k_N \to\infty$ and $p_N= o(N)$. Let $y_{i,j}^{(N)} = \sqrt{N} (X_N)_{i,j}$ be the scaled entries for $1\le i \le m_N, 1\le j \le k_N$. Define the random empirical measure by
		$$
		\nu_N := \frac{1}{p_N} \sum_{i=1}^{m_N} \sum_{j=1}^{k_N} \delta_{y_{i,j}^{(N)}}.
		$$
		Then the sequence $\nu_N$ satisfies an LDP on the space of probability measures $\mathcal{P}(\R)$, equipped with the weak topology, with speed $p_N = m_N k_N$ and good convex rate function $I$ given by the relative entropy with respect to the standard Gaussian measure $\gamma$:
		$$
		I(\mu) = H(\mu | \gamma) := \begin{cases} \int_{\R} \log \left( \frac{d\mu}{d\gamma} \right) d \mu & : \mu \ll \gamma, \\ +\infty & :\text{otherwise.} \end{cases}
		$$
	  \end{thmalpha}

	  The proof of the theorem, as well as the proofs of most other results below, is given in Section \ref{sec:proofs}. For example, the LDP upper bound gives
	  $$
	  \limsup_{N\to\infty}\frac{1}{p_N}\log\P\left(\frac{1}{p_N} \sum_{i=1}^{m_N} \sum_{j=1}^{k_N} \delta_{y_{i,j}^{(N)}} \not\in \mathbb{B}_\epsilon(\gamma)\right)
	  \le -\inf_{\mu\notin \mathbb{B}_\epsilon(\gamma)}H(\mu|\gamma)
	  $$
	  for $\epsilon > 0$, where $\mathbb{B}_\epsilon(\gamma)$ is the ball of radius $\epsilon$ around $\gamma$ in the Lévy-Prokhorov metric.

	  \begin{rmk}
		The condition $p_N = o(N)$ cannot be improved without changing the rate function. As we have seen in the introduction, by \cite[Theorem 2.8]{kim2021large}, $\frac{1}{N} \sum_{i=1}^{N} \delta_{\sqrt{N} v_i}$, where $v_i$ are the columns of $V_N = (v_1, \ldots, v_N) \sim \Unif(\V_{m,N})$, satisfies an LDP in the $q$-Wasserstein topology for $q\in (0,2)$ (and hence in the weak topology) with speed $N$ and good rate function given by Equation \eqref{eq:ratefunctionstiefelkim}. Using the contraction principle (see Proposition \ref{prop:contractionprinciple_main} below), for $p_N = kN$ ($k$ fixed), the sequence $\nu_N$ satisfies an LDP with speed $N$ and a rate function that is not just given by the relative entropy but also involves the covariance operator of $\nu_N$.
	  \end{rmk}

	  \noindent The theorem leads to the following corollary:

	  \begin{cor}
		Assume that the matrices in Theorem \ref{thm:empirical_main} are defined on a common probability space and that $p_N=\lfloor N^\alpha\rfloor$ for some $\alpha \in (0,1)$. Then, with probability $1$, the sequence of empirical measures $\nu_N$ converges to the standard Gaussian measure $\gamma$ as $N\to\infty$ in the weak topology.
	  \end{cor}
	  \begin{proof}
		Let $B$ be an open neighborhood of $\gamma$ in $\mathcal{P}(\R)$. Since $H(\,\cdot\,|\gamma)$ is a good rate function with unique zero $\gamma$, we have $c_B:=\inf_{\mu\in B^c}H(\mu|\gamma)>0$. Indeed, otherwise there would be a sequence $(\mu_n)_{n\in\N}$ in $B^c$ such that $H(\mu_n|\gamma)\to0$. The compactness of the sublevel sets and lower semicontinuity would give a limit point in $B^c$ at which the relative entropy vanishes, contradicting $\gamma\in B$. The LDP upper bound implies that, for all sufficiently large $N$,
		$$
		\P(\nu_N\notin B)\le \exp\left(-\frac{c_Bp_N}{2}\right).
		$$
		Hence $\sum_{N\in\N}\P(\nu_N\notin B)<\infty$. The Borel-Cantelli lemma shows that $\nu_N$ is eventually in $B$ with probability $1$. Applying this argument to the countable neighborhood basis $(\mathbb{B}_{1/r}(\gamma))_{r\in\N}$ and intersecting the resulting events of probability $1$ proves the claim.
	  \end{proof}

	  \noindent The proof of Theorem \ref{thm:empirical_main} uses the following result on the uniform convergence of density functions which might be of independent interest. Let \(m=m_N\), \(k=k_N\) and \(N\ge m+k\). Define
		\[
		Z_N := \bigl(Z_{ij}^{(N)}\bigr)_{1\le i\le m,\;1\le j\le k}
		\]
		to be the upper-left \(m\times k\) block of a Haar-distributed matrix in \(\O(N)\). \noindent By \cite[Lemma 2.5]{jiang2006many}, the block $Z_N$ has, for $N\ge k+m$, the density $f_N:[-1,1]^{m\times k}\to [0,\infty)$,
		$$
		f_N(A)=\frac{\Gamma_m\left(\frac{N}{2}\right)}{\pi^{\frac{m k}{2}} \Gamma_m\left(\frac{N-k}{2}\right)} \operatorname{det}\left(\operatorname{Id}_{m}-A A^\top\right)^{\frac{N-k-m-1}{2}} \mathbf{1}_{\left\{\left\|A A^\top\right\|_{op}<1\right\}}, \quad A \in[-1,1]^{m \times k},
		$$
		where $\Gamma_m$ denotes the multivariate gamma function
		$$
		\Gamma_m(x)=\pi^{\frac{m(m-1)}{4}} \prod_{i=1}^m \Gamma\left(x-\frac{i-1}{2}\right),
		$$
		where $x> \frac{1}{2}(m-1)$ (see \cite[Theorem 1.4.1]{gupta2018matrix}). Here \(\|AA^\top\|_{op}\) is the operator norm of \(AA^\top\). For \(B\in\mathbb R^{m\times k}\) define
		\[
		g_N(B)
		:= N^{-\frac{mk}{2}}\,f_N\left(\frac{B}{\sqrt{N}}\right),
		\qquad
		\phi_N(B)
		:= (2\pi)^{-\frac{mk}{2}}
		\exp\Bigl(-\frac12\|B\|_F^2\Bigr),
		\]
		where $\|\cdot\|_F$ denotes the Frobenius norm. Note that $g_N$ is the density of $\sqrt{N} Z_N$ while $\phi_N$ is the density of the standard Gaussian $m\times k$ block. We prove the following result:

		\begin{proposition}[Uniform Local Limit]
		\label{prop:uniform_local_limit}
		Assume \(m, k \ge 1\) with \(m = m_N = o(N)\) and \(k= k_N = o(N)\). Let \(R_N\) be a sequence such that \(R_N = o(N^{1/2})\). Then there exists a constant \(C>0\) such that for all sufficiently large \(N\),
		\[
		\sup_{B\in\mathbb R^{m\times k}: \|B\|_F \le R_N}
		\Bigl|
		\log \frac{g_N(B)}{\phi_N(B)}
		\Bigr|
		\;\le\;
		C\Bigl(\frac{R_N^4}{N} + \frac{mk(m+k)}{N} + \frac{(k+m)R_N^2}{N}\Bigr).
		\]
		\end{proposition}

\begin{rmk}
Under stronger conditions, specifically if $mk(m+k)=o(N)$ and $R_N = o(N^{1/4})$, the error term on the right hand side tends to zero as $N \to \infty$. This implies that $g_N(\cdot)/\phi_N(\cdot) \to 1$ uniformly on $\{B \in \mathbb{R}^{m\times k} : \|B\|_F \le R_N\}$.
\end{rmk}

    The next main results are the MDPs for scaled orthogonal and Stiefel matrices, which complement the LDP results in \cite[Theorem A and B]{KabluchkoLarge}.

	\begin{thmalpha}[MDP for the Stiefel Manifold]\label{thm:MDPforStiefel_main} 
		Let $m\in\N$ be fixed and $V_{m,N}$ be chosen uniformly (i.e. with respect to the Haar measure) in $\mathbb{V}_{m,N}$ for every $N\ge m$. Embed $V_{m,N}$ in $\R^{m\times\infty}$ by padding it with zero columns, and equip this space with the product topology. Furthermore, let $\beta_N$ be a sequence of positive numbers such that $\beta_N\to\infty$ and $\beta_N = o(\sqrt{N})$. Then, the sequence $\beta_N V_{m,N}$ satisfies an LDP on the space $\R^{m\times\infty}$ with speed $\frac{N}{\beta_N^2}$ and good convex rate function 
		$$
		I:\R^{m\times\infty}\to [0,\infty],\quad I(A) = \frac{1}{2}\|A\|_F^2.
		$$ 
		\end{thmalpha}

		A similar result is true for the orthogonal group.
	
		\begin{thmalpha}[MDP for the Orthogonal Group]\label{thm:MDPorthogonalgroup} 
		Let $A_{N}$ be chosen uniformly in $\O(N)$ for every $N$. Embed $A_N$ in $\R^{\infty\times\infty}$ by padding it with zero rows and columns, and equip this space with the product topology. Furthermore, let $\beta_N$ be a sequence of positive numbers such that $\beta_N\to\infty$ and $\beta_N = o(\sqrt{N})$. Then, the sequence $\beta_N A_{N}$ satisfies an LDP on the space $\R^{\infty\times\infty}$ with speed $\frac{N}{\beta_N^2}$ and good convex rate function 
		$$
		J:\R^{\infty\times\infty}\to [0,\infty],\quad J(A) = \frac{1}{2}\|A\|_F^2.
		$$ 
		\end{thmalpha}
	
	\noindent These theorems show the asymptotic behavior below the CLT level.  Deviations from the zero matrix are governed by a Gaussian-like quadratic rate function and speed $\frac{N}{\beta_N^2}$ which interpolates between the LDP speed $N$ and the CLT.


\section{Preliminaries}\label{sec:preliminaries}\label{C:sec:preliminaries}

We introduce some notation. The space of probability measures in $\R$ is denoted by $\mathcal{P}(\R)$. We write $\mathcal{N}(0,\operatorname{Id}_m)$ for the standard $m$-dimensional Gaussian distribution on $\R^m$, where $\operatorname{Id}_m$ is the $m$-dimensional identity matrix. We denote the Dirac measure at some point $x\in \R^m$ by $\delta_x$. For a random vector $X$ that is distributed according to a probability measure $\mu$, we write $X\sim \mu$. The mean of $X$ is denoted by $\E[X]$ and the variance by $\Var[X]$. An open ball of radius $r > 0$ centered at $x$ in a metric space is denoted by $\mathbb{B}_r(x)$. The \emph{weak topology} is the coarsest topology on $\mathcal{P}(\R)$ for which the maps $\mu \mapsto \int f d\mu$ are continuous for all $f \in C_b(\R)$, the space of bounded continuous functions. As mentioned above, $\|\cdot\|_F$ and $\|\cdot\|_{op}$ are the Frobenius norm and operator norm, respectively. Depending on the context, we use $\mathbf{1}_{A}(x)$ or $\mathbf{1}\{x\in A\}$ for the indicator function over some set $A$. For sequences $a_N$, $b_N$, we write $a_N = o(b_N)$ if $|a_N/b_N| \to 0$ and $a_N = O(b_N)$ if there exist $C>0$ and $N_0\in\N$ such that $|a_N|\le C|b_N|$ for all $N\ge N_0$. We use $a_N \sim b_N$ to denote asymptotic equivalence, i.e. $\frac{a_N}{b_N}\to 1$. For a square matrix $A$, $\Tr{A}$ denotes the trace of $A$ and for a measurable set $B$, $\vol_m(B)$ is the $m$-dimensional Lebesgue measure of $B$. 
	
	\subsection{Large and Moderate Deviation Principles}

	In this section, we collect some preliminary definitions and results that are used in the proofs of the main results. We start with the concept of a large deviation principle. Let $\mathbb{X}$ be a topological space. A function $I:\mathbb{X}\to[0,\infty]$ is a \emph{rate function} if it is lower semi-continuous, i.e. its sublevel sets $\{x\in\mathbb{X}: I(x)\le y\}$ are closed for all $y \ge 0$. If those sets are, in addition, compact, we call $I$ a \emph{good rate function}.

	\begin{df}[Large Deviation Principle]\label{def:LDP_main} 
		A sequence of random variables $X_N$ taking values in a topological space $\mathbb{X}$ satisfies a \emph{Large Deviation Principle (LDP)} with \emph{speed} $s_N \to \infty$ and \emph{rate function} $I:\mathbb{X}\to [0,\infty]$ if $I$ is a rate function and
		\begin{align}
			\liminf_{N\to\infty} \frac{1}{s_N}\log \P(X_N \in O) &\ge -\inf_{x\in O}I(x) \quad \text{ for all open } O\subseteq\mathbb{X},\label{LDPlowerbound_main}\\ 
			\limsup_{N\to\infty} \frac{1}{s_N}\log \P(X_N\in C) &\le -\inf_{x\in C}I(x) \quad \text{ for all closed } C\subseteq\mathbb{X}\label{LDPupperbound_main}.
		\end{align}
	\end{df}

	We say that a sequence satisfies a \emph{weak LDP} with speed $s_N$ and rate function $I$ if it satisfies the assumptions of an LDP with speed $s_N$ and rate function $I$ except that the upper bound in \eqref{LDPupperbound_main} is only assumed to hold for compact sets. Next, we introduce the notion of a moderate deviation principle which formally is just a special kind of LDP. 

    \begin{df}[Moderate Deviation Principle (MDP)]\label{def:MDP_main} 
        Let $Y_N$ be a sequence of random variables. Let $\beta_N$ be a sequence such that $\beta_N \to \infty$ and $\beta_N = o(\sqrt{N})$. We say that $\beta_N Y_N$ satisfies a \emph{Moderate Deviation Principle (MDP)} with speed $\frac{N}{\beta_N^2}$ and rate function $I$ if it obeys an LDP with the same speed and rate function.
    \end{df}

	\noindent Since the speed $\frac{N}{\beta_N^2}$ tends to infinity and is $o(N)$, the MDP scale lies between the central limit theorem (CLT) scale and the LDP at speed $N$. In the MDP case, the sequence still satisfies an LDP but the CLT behavior is often reflected in the rate function. We can see this in the Theorems \ref{thm:MDPforStiefel_main} and \ref{thm:MDPorthogonalgroup} as the rate function has a quadratic form that reminds of the exponent in the Gaussian density. 

	\begin{df}[Exponential Tightness]\label{def:exponentialtightness_main} 
		A sequence of random variables $X_N$ with values in a topological space $\mathbb{X}$ is \emph{exponentially tight} with respect to speed $s_N$ if for every $a>0$, there exists a compact set $K\subseteq \mathbb{X}$ such that 
		$$\limsup_{N\to\infty} \frac{1}{s_N} \log \P(X_N\notin K) \le -a.$$
	\end{df}

	\noindent We continue with some standard tools in LDP theory. The following result is taken from Lemma 1.2.18 in \cite{dembo2009techniques} and the remark before.
	\begin{proposition}[Equivalence of LDP and weak LDP \& exponential tightness]\label{thm:LDPequivalence} 
		Let $\mathbb{X}$ be a locally compact space or a Polish space. Then a sequence of random variables in $\mathbb{X}$ satisfies an LDP with speed $s_N$ and good rate function $I$ if and only if it satisfies the corresponding weak LDP and is exponentially tight. 
	\end{proposition}

\noindent The next proposition is adapted from \cite[Lemma 4.1.11 and Lemma 4.1.18]{dembo2009techniques}.
\begin{proposition}[Reformulation of weak LDP]\label{prop:weak-ldp-base}
		Let $\mathcal{T}$ be a base of the topology in a metric space $\mathbb{X}$. Let $X_N$ be a sequence of random variables with values in $\mathbb{X}$ and assume $s_N$ is some sequence with $s_N\to\infty$. If for every $x\in \mathbb{X}$,
		\[
		  I(x)
		  := - \inf_{A\in\mathcal{T} : x\in A}
			   \limsup_{N\to\infty}\frac{1}{s_N}\log \mathbb{P}(X_N\in A)
			= - \inf_{A\in\mathcal{T} : x\in A}
			   \liminf_{N\to\infty}\frac{1}{s_N}\log \mathbb{P}(X_N\in A),
		\]
		then $X_N$ satisfies a weak LDP with speed $s_N$ and rate function $I$.  Conversely, if $X_N$ satisfies a weak LDP with speed $s_N$ and rate function $I$, then the above identities hold.
	\end{proposition}

\noindent The next proposition, known as the Contraction Principle, can be found in \cite[Theorem 4.2.1]{dembo2009techniques}.

\begin{proposition}[Contraction Principle]\label{prop:contractionprinciple_main}
			Let $f:\mathbb{X}\to\mathbb{Y}$ be continuous, where $\mathbb{X}$ and $\mathbb{Y}$ are Hausdorff spaces. If $X_{N}$ satisfies an LDP on $\mathbb{X}$ with speed $s_N$ and good rate function $I$, then $f(X_N)$ satisfies an LDP on $\mathbb{Y}$ with speed $s_N$ and good rate function $J(y) = \inf_{x\in f^{-1}(y)}I(x)$ (with the convention that $\inf_{x\in \emptyset}I(x) = +\infty$). 
	\end{proposition}

	We introduce the notion of projective limits which we will need for the proofs of Theorems \ref{thm:MDPforStiefel_main} and \ref{thm:MDPorthogonalgroup} (projective limits were also used in the proof of \cite[Theorem A]{KabluchkoLarge}). A \emph{projective system} $\left(\mathbb{Y}_j, p_{i j}\right)_{i \leq j}$ consists of Hausdorff topological spaces $\mathbb{Y}_j, j \in \mathbb{N}$, and continuous mappings $p_{i j}: \mathbb{Y}_j \rightarrow \mathbb{Y}_i$ such that for all $i \leq j \leq k$, we have $p_{i k}=p_{i j} \circ p_{j k}$ and $p_{j j}, j \in \mathbb{N}$ is the identity mapping on $\mathbb{Y}_j$. Then, the projective limit $\mathbb{Y}$ of this system is given by
	$$
\mathbb{Y}:=\Big\{y=\left(y_j\right)_{j \in \mathbb{N}} \in \prod_{j \in \mathbb{N}} \mathbb{Y}_j \,:\, y_i=p_{i j}\left(y_j\right), \forall i<j\Big\}.
        $$
        We denote by $p_j:\mathbb{Y}\to \mathbb{Y}_j$ for $j\in\N$ the canonical projection, i.e., for $y\in\mathbb{Y}$, we have $p_j(y) = y_j$ in the notation above. The following theorem by Dawson and G\"artner, see \cite[Theorem 3.3]{dawsont1987large}, allows inferring an LDP on $\mathbb{Y}$ given an LDP on $\mathbb{Y}_j$ for each $j$.
	
\begin{proposition}[Dawson -- G\"artner]\label{prop:B:dawson-gaertner}
		Let $\mathbb{Y}$ be the projective limit of the projective system $\left(\mathbb{Y}_j, p_{i j}\right)_{i \leq j}$. Assume that $X_N$ is a sequence of random variables in $\mathbb{Y}$ such that for any $\ell \in \mathbb{N}$, $p_{\ell}(X_N)$ satisfies an LDP on $\mathbb{Y}_{\ell}$ with speed $s_N$ and good rate function $I_{\ell}$. Then $X_N$ satisfies an LDP with speed $s_N$ and good rate function $I: \mathbb{Y} \rightarrow[0,+\infty]$ given by
		$$
		I(x):=\sup _{\ell \in \mathbb{N}} I_{\ell}\left(p_{\ell}(x)\right) .
		$$
	\end{proposition}

\section{Proofs}
\label{C:sec:proofs}\label{sec:proofs}

In this section, we provide the proofs, starting with the proof of the LDP for the empirical measure of scaled entries of an orthogonal matrix, followed by the proofs of the MDPs for scaled orthogonal and Stiefel matrices.

It is possible to prove Theorem \ref{thm:empirical_main} by establishing exponential equivalence with the Gaussian empirical measures. Such a proof requires steps similar to those below and does not seem to reduce the lenght substantially. We use the density approach because the uniform estimate in Proposition \ref{prop:uniform_local_limit} may also be useful in other settings.

	\subsection{Proof of Theorem \ref{thm:empirical_main}}
	We start with the following lemma that will be essential in the proof of Proposition \ref{prop:uniform_local_limit}. It quantifies the difference in the exponents of both densities. Recall that we compare the density $g_N(B)=N^{-mk/2}f_N(B/\sqrt{N})$ with
\[
  f_N(A)=\frac{\Gamma_m\left(\frac{N}{2}\right)}{\pi^{\frac{mk}{2}}\,\Gamma_m\left(\frac{N-k}{2}\right)}\,\det\bigl(\operatorname{Id}_m-AA^T\bigr)^{\frac{N-k-m-1}{2}}\;\mathbf 1_{\{\|AA^T\|_{op}<1\}},
\]
	with the Gaussian density $\phi_N(B)=(2\pi)^{-mk/2}\exp(-\frac12\|B\|_F^2)$.

\begin{lemma}[Log-determinant expansion]
	\label{lem:log_det_expansion}
	Let \(X := \frac1N\,B B^\top\) where $B$ is a real $m\times k$ matrix. For \(\|X\|_{op}<1\), the following expansion holds:
	\[
	\frac{N-k-m-1}{2}\,\log\det(\operatorname{Id}_m - X)
	= -\frac{\|B\|_F^2}{2}
	  + \frac{k+m+1}{2N}\,\|B\|_F^2
	  - \frac{\mathrm{Tr}\bigl((BB^\top)^2\bigr)}{4N}
	  + R_N'(B).
	\]
	If \(\|B\|_F \le R_N\) for a sequence $R_N$ with \(R_N = o(N^{1/2})\), the remainder term is bounded by
	\[
	|R_N'(B)| \;\le\; C \left( \frac{(k+m)R_N^4}{N^2} + \frac{R_N^6}{N^2} \right)
	\]
	for some constant $C>0$.
\end{lemma}

\begin{proof}
	We start with an identity for $\log\det(\operatorname{Id}_m - X)$. By \cite[Equation (512)]{petersen2008matrix}, we have 
  $$
  \det(\operatorname{Id}_m - X) = \exp(\mathrm{Tr}(\log(\operatorname{Id}_m -X)))
  $$ 
  and thus $\log(\det(\operatorname{Id}_m - X)) = \mathrm{Tr}(\log(\operatorname{Id}_m - X))$. Since $\|X\|_{op}<1$, we get by a series expansion of the logarithm and linearity of $\mathrm{Tr}$,
	$$ \log\det(\operatorname{Id}_m - X) = -\sum_{r=1}^\infty\frac{\mathrm{Tr}(X^r)}{r}. $$
	Let \(C_N = \frac{N-k-m-1}{2}\) and note that $\mathrm{Tr}(X) = \|B\|_F^2/N$. Then
	\begin{align*} 
    &C_N \log\det(\operatorname{Id}_m-X) = C_N \left(-\mathrm{Tr}(X) - \frac12\mathrm{Tr}(X^2) - \sum_{r\ge3}\frac{\mathrm{Tr}(X^r)}{r}\right) \\
	&= \frac{N}{2}(1-\frac{k+m+1}{N}) \left(-\frac{1}{N}\|B\|_F^2 - \frac{1}{2N^2}\mathrm{Tr}((BB^\top)^2) - \sum_{r\ge3}\frac{\mathrm{Tr}(X^r)}{r}\right) \\
	&= -\frac12\|B\|_F^2 - \frac{1}{4N}\mathrm{Tr}((BB^\top)^2) + \frac{k+m+1}{2N}\|B\|_F^2 + \frac{k+m+1}{4N^2}\mathrm{Tr}((BB^\top)^2)- C_N \sum_{r\ge3}\frac{\mathrm{Tr}(X^r)}{r}.
	\end{align*}
	The remainder term \(R_N'(B)\) is therefore given by
	\[ R_N'(B) = \frac{k+m+1}{4N^2}\mathrm{Tr}((BB^\top)^2) - C_N \sum_{r\ge3}\frac{\mathrm{Tr}(X^r)}{r}. \]
	We bound the two terms. Let us start with the first term, \(\frac{k+m+1}{4N^2}\mathrm{Tr}((BB^\top)^2)\). We use the inequality \(\mathrm{Tr}((BB^\top)^2) \le (\|B\|_F^2)^2 = \|B\|_F^4\). To see this, let \(\sigma_i \ge 0\) be the eigenvalues of the \(m \times m\) matrix \(BB^\top\). Then \(\mathrm{Tr}(BB^\top) = \sum_{i=1}^m \sigma_i = \|B\|_F^2\), and \(\mathrm{Tr}((BB^\top)^2) = \sum \sigma_i^2\le (\sum \sigma_i)^2\). Thus, this first term of \(R_N'(B)\) is bounded by:
	$$ D_1 \frac{(k+m)\|B\|_F^4}{N^2}, $$
	where $D_1>0$ is some constant. Next, we analyze the second term. We have $\|X\|_{op} = \frac{1}{N}\|BB^\top\|_{op} \le \frac{1}{N}\|B\|_{op}^2 \le \frac{1}{N}\|B\|_F^2$ by standard properties of the operator norm. Since $\|B\|_F \le R_N = o(\sqrt{N})$, we have $\|X\|_{op} = o(1)$. Next, for the sum term \(-C_N \sum_{r\ge3}\frac{\mathrm{Tr}(X^r)}{r}\), we use the inequality \(\mathrm{Tr}(X^r) \le \|X\|_{op}^{r-2}\mathrm{Tr}(X^2)\) for \(r \ge 2\). This holds because $\|X\|_{op} = \lambda_1$ and for \(r \ge 2\), \(\lambda_i^r = \lambda_i^{r-2}\lambda_i^2 \le \lambda_1^{r-2}\lambda_i^2 = \|X\|_{op}^{r-2}\lambda_i^2\) where $\lambda_1\ge \ldots\ge\lambda_m \ge 0$ are the eigenvalues of $X$. Summing over \(i\) gives \(\mathrm{Tr}(X^r) = \sum_i \lambda_i^r \le \|X\|_{op}^{r-2} \sum_i \lambda_i^2 = \|X\|_{op}^{r-2}\mathrm{Tr}(X^2)\). Thus,
	$$\Big|C_N \sum_{r\ge3}\frac{\mathrm{Tr}(X^r)}{r}\Big| \le C_N \sum_{r\ge3}\frac{|\mathrm{Tr}(X^r)|}{r} \le C_N \mathrm{Tr}(X^2) \sum_{r\ge3}\frac{\|X\|_{op}^{r-2}}{r}.$$
	Using the geometric series and the fact that \(\|X\|_{op} =o(1)\), we get \(\sum_{r\ge3}\frac{\|X\|_{op}^{r-2}}{r} \le \frac{1}{3} \frac{\|X\|_{op}}{1-\|X\|_{op}}\).
	This leads to the overall bound for the sum term:
	$$ \Big|C_N \sum_{r\ge3}\frac{\mathrm{Tr}(X^r)}{r}\Big| \le C_N \mathrm{Tr}(X^2) \frac{\|X\|_{op}}{3(1-\|X\|_{op})}. $$
	We have established above that \(\|X\|_{op} \le N^{-1}\|B\|_F^2\) and \(\mathrm{Tr}(X^2) = N^{-2}\mathrm{Tr}((BB^\top)^2) \le N^{-2}\|B\|_F^4\). Since \(\|X\|_{op}=o(1)\), we have \(1/(3(1-\|X\|_{op})) = O(1)\). Substituting these approximations, we get
	$$ C_N \mathrm{Tr}(X^2) \frac{\|X\|_{op}}{3(1-\|X\|_{op})} =  O\left(\frac{\|B\|_F^6}{N^2}\right). $$
	Combining these gives
	$$ |R_N'(B)| \le C \left(\frac{(k+m)\|B\|_F^4}{N^2} + \frac{\|B\|_F^6}{N^2}\right). $$
	Since \(\|B\|_F \le R_N\), we get the desired bound
	$$ |R_N'(B)| \le C \left(\frac{(k+m)R_N^4}{N^2} + \frac{R_N^6}{N^2}\right). $$
	\end{proof}

	We continue with a technical lemma about the asymptotics of the multivariate gamma function. This compares the normalizing constants of both densities. 

	\begin{lemma}[Gamma quotient estimate]
		\label{lem:gamma_quotient}
		Let \(\Gamma_m(x) = \pi^{\frac{m(m-1)}4}\prod_{j=1}^m\Gamma\bigl(x-\frac{j-1}2\bigr)\) for $x > (m-1)/2$ be the multivariate gamma function. Assume \(k = k_N = o(N)\) and \(m = m_N = o(N)\). Then
		\[
		\log\frac{\Gamma_m(\frac N2)}{\Gamma_m(\frac{N-k}2)}
		  - \frac{mk}{2}\,\log\Bigl(\frac N2\Bigr)
		= O\Bigl(\frac{mk(m+k)}{N}\Bigr).
		\]
		\end{lemma}

		\begin{proof}
			Set $x_l=(N-k-l)/2$ for $l=0,\ldots,m-1$ and $h=k/2$. The assumptions imply that $x_l+t$ is of order $N$ uniformly for $0\le l\le m-1$ and $0\le t\le h$. Using
			\[
			\log\Gamma(x_l+h)-\log\Gamma(x_l)=\int_0^h\psi(x_l+t)\,\dint t
			\]
			and the estimate $\psi(x)=\log x+O(x^{-1})$ as $x\to\infty$ from \cite[Section 5.11]{lozier2003nist}, we obtain, uniformly in $l$,
			\[
			\log\Gamma(x_l+h)-\log\Gamma(x_l)
			=\int_0^h\log(x_l+t)\,\dint t+O\left(\frac{h}{N}\right).
			\]
			For $0\le t\le h$,
			\[
			\log(x_l+t)
			=\log\left(\frac{N}{2}\right)
			+\log\left(1-\frac{k+l-2t}{N}\right)
			=\log\left(\frac{N}{2}\right)+O\left(\frac{k+l}{N}\right)
			\]
			uniformly in $l$ and $t$. Consequently,
			\[
			\log\Gamma(x_l+h)-\log\Gamma(x_l)
			=\frac{k}{2}\log\left(\frac{N}{2}\right)
			+O\left(\frac{k(k+l+1)}{N}\right).
			\]
			Summing over $l=0,\ldots,m-1$ gives
			\[
			\log\frac{\Gamma_m(\frac{N}{2})}{\Gamma_m(\frac{N-k}{2})}
			=\frac{mk}{2}\log\left(\frac{N}{2}\right)
			+O\left(\frac{mk(m+k)}{N}\right),
			\]
			which proves the claim.
		\end{proof}

		We now combine both lemmas to prove Proposition \ref{prop:uniform_local_limit}.

	\begin{proof}[Proof of Proposition~\ref{prop:uniform_local_limit}]
		First, we confirm that the indicator is $1$ for large $N$, i.e. \(\|AA^\top \|_{op}<1\) holds for \(\|B\|_F \le R_N\). With $B = \sqrt{N}A$ we have \(AA^\top = \frac{1}{N}BB^\top\) and thus \(\|AA^\top\|_{op} = \frac{1}{N}\|BB^\top\|_{op}\).  Using \(\|BB^\top\|_{op} \le \|B\|_F^2\) and the assumption \(\|B\|_F \le R_N\), we have
		\[ \|AA^\top\|_{op} \le \frac{1}{N} \|B\|_F^2 \le \frac{R_N^2}{N} = o(1). \]
		Hence \(\mathbf1_{\{\|AA^\top\|_{op}<1\}} = 1\) for large enough \(N\). Using logarithm rules, 
		\begin{align*}
		 \log g_N(B) &= -\frac{mk}{2}\log N + \log\Gamma_m\bigl(\frac N2\bigr) - \log\Gamma_m\bigl(\frac{N-k}2\bigr) - \frac{mk}{2}\log\pi\\
		 & + \frac{N-k-m-1}2\log\det(\operatorname{Id}_m - \frac{BB^\top}N), 
		\end{align*}
		and
		$$ \log\phi_{N}(B) = -\frac{mk}{2}\log(2\pi) - \frac12\,\|B\|_F^2 = -\frac{mk}{2}\log 2 - \frac{mk}{2}\log\pi - \frac12\,\|B\|_F^2. $$
		Hence, with \(X=\frac1N BB^\top\),
		\begin{align*} \log\frac{g_N(B)}{\phi_{N}(B)} &= \left[ \log\frac{\Gamma_m(\frac N2)}{\Gamma_m(\frac{N-k}2)} - \frac{mk}{2}\log N + \frac{mk}{2}\log 2 \right]  \\ & \quad + \left[ \frac{N-k-m-1}{2}\log\det(\operatorname{Id}_m-X) + \frac12\|B\|_F^2 \right].
		\end{align*}
		By Lemma~\ref{lem:gamma_quotient},
		$$ \log\frac{\Gamma_m(\frac N2)}{\Gamma_m(\frac{N-k}2)} = \frac{mk}{2}\log(\frac N2) + O\left(\frac{mk(m+k)}{N}\right), $$
		so the first term in brackets becomes
		$$ \frac{mk}{2}\log(\frac N2) + O\left(\frac{mk(m+k)}{N}\right) - \frac{mk}{2}\log N + \frac{mk}{2}\log 2 = O\left(\frac{mk(m+k)}{N}\right). $$
		By Lemma~\ref{lem:log_det_expansion}, the second term is
		$$ \frac{k+m+1}{2N}\|B\|_F^2 - \frac{\mathrm{Tr}[(BB^\top)^2]}{4N} + R_N'(B). $$
		Therefore, for \(\|B\|_F\le R_N\), using \(\mathrm{Tr}[(BB^\top)^2] \le \|B\|_F^4 \le R_N^4\) again, for large $N$,
		\begin{align*} 
		\Bigl|\log\frac{g_N(B)}{\phi_{N}(B)}\Bigr| &\le C_1\frac{mk(m+k)}{N} + \frac{k+m+1}{2N}R_N^2 + \frac{R_N^4}{4N} + |R_N'(B)|\\
		& \le C_1\frac{mk(m+k)}{N} + C_2\frac{(k+m)R_N^2}{N} + C_3\frac{R_N^4}{N} + C_4\left(\frac{(k+m)R_N^4}{N^2} + \frac{R_N^6}{N^2}\right)\\
		& \le C_5\left(\frac{R_N^4}{N} + \frac{mk(m+k)}{N} + \frac{(k+m)R_N^2}{N}\right)
		\end{align*}
		for positive constants \(C_1, C_2, C_3, C_4, C_5\). The last inequality uses $(m+k)/N=o(1)$ and $R_N^2/N=o(1)$. 
		\end{proof}

		We now compare probabilities under the Haar density and the Gaussian density on a growing Frobenius ball.

		\begin{lemma}\label{lem:density_comparison}
			With the notation of Theorem \ref{thm:empirical_main}, let $Y_N=(y_{i,j}^{(N)})_{i,j}$ and let $G_N=(g_{i,j})_{i,j}$ be an $m_N\times k_N$ matrix with i.i.d. standard Gaussian entries. Define $L_N \colon \R^{m_N\times k_N}\to\mathcal P(\R)$ by
			$$
			L_N(B)=\frac{1}{p_N}\sum_{i=1}^{m_N}\sum_{j=1}^{k_N}\delta_{b_{i,j}}.
			$$
			There exist sequences $R_N>0$ and $\eta_N\to0$ such that, for every measurable set $A\subseteq\mathcal P(\R)$,
			\begin{align*}
			\P(L_N(Y_N)\in A)&\le e^{\eta_Np_N}\P(L_N(G_N)\in A)+\P(\|Y_N\|_F>R_N),\\
			\P(L_N(Y_N)\in A)&\ge e^{-\eta_Np_N}\bigl(\P(L_N(G_N)\in A)-\P(\|G_N\|_F>R_N)\bigr).
			\end{align*}
			Moreover,
			$$
			\limsup_{N\to\infty}\frac{1}{p_N}\log\P(\|Y_N\|_F>R_N)
			=
			\limsup_{N\to\infty}\frac{1}{p_N}\log\P(\|G_N\|_F>R_N)
			=-\infty.
			$$
		\end{lemma}

		\begin{proof}
			Set $R_N=N^{1/8}p_N^{3/8}$. Then $R_N=o(N^{1/2})$, $R_N^2/p_N\to\infty$ and $R_N^4/(Np_N)\to0$. Proposition \ref{prop:uniform_local_limit} and $m_N+k_N\le p_N+1$ imply that
			$$
			\eta_N:=\frac{1}{p_N}\sup_{\|B\|_F\le R_N}\left|\log\frac{g_N(B)}{\phi_N(B)}\right|\longrightarrow0.
			$$
			Indeed, after division by $p_N$, the three terms in Proposition \ref{prop:uniform_local_limit} are bounded by constant multiples of $(p_N/N)^{1/2}$, $(m_N+k_N)/N$ and
			$$
			\frac{(m_N+k_N)R_N^2}{Np_N}\le \left(1+\frac{1}{p_N}\right)\left(\frac{p_N}{N}\right)^{3/4},
			$$
			respectively. Each expression tends to zero.
			The moment generating function of the chi-square distribution gives
			$$
			\P(\|G_N\|_F>R_N)\le \exp\left(-\frac{R_N^2}{4}+\frac{p_N\log2}{2}\right).
			$$
			Also, $\E[\|Y_N\|_F]\le\sqrt{p_N}$. The function $F_N \colon \O(N)\to\R$, $F_N(U)=\|\sqrt{N}U^{m_N,k_N}\|_F$, is $\sqrt{N}$ Lipschitz with respect to the Hilbert-Schmidt distance. The two components of $\O(N)$ give the same distribution for $F_N$, since multiplication by a diagonal reflection in a row outside the block changes the determinant and leaves $F_N$ unchanged. Since $R_N/\sqrt{p_N}\to\infty$, we have $R_N-\E[F_N]\ge R_N/2$ for all sufficiently large $N$. Proposition 10 in \cite{meckes2017rates} applied on either component therefore yields $\P(\|Y_N\|_F>R_N)\le e^{-cR_N^2}$ for all sufficiently large $N$ and some $c>0$. Thus, both tail probabilities are negligible at speed $p_N$ in the sense of the last display in the statement.

			On the ball $\{B:\|B\|_F\le R_N\}$, the definition of $\eta_N$ gives
			$$
			e^{-\eta_Np_N}\phi_N(B)\le g_N(B)\le e^{\eta_Np_N}\phi_N(B).
			$$
			Integrating these inequalities over $\{B:L_N(B)\in A,\ \|B\|_F\le R_N\}$ gives the two asserted bounds.
		\end{proof}

\begin{proof}[Proof of Theorem \ref{thm:empirical_main}]
	Let $\widehat\nu_N=L_N(G_N)$. By Sanov's theorem \cite[Theorem 6.2.10]{dembo2009techniques}, the sequence $(\widehat\nu_N)_{N\in\N}$ satisfies an LDP with speed $p_N$ and good convex rate function $H(\,\cdot\,|\gamma)$.

	Let $C\subseteq\mathcal P(\R)$ be closed. The first inequality in Lemma \ref{lem:density_comparison} and the elementary logarithmic estimate for a sum give
	\begin{align*}
	\limsup_{N\to\infty}\frac{1}{p_N}\log\P(\nu_N\in C)
	&\le \max\left\{\limsup_{N\to\infty}\left(\eta_N+\frac{1}{p_N}\log\P(\widehat\nu_N\in C)\right),
	\limsup_{N\to\infty}\frac{1}{p_N}\log\P(\|Y_N\|_F>R_N)\right\}\\
	&\le \max\left\{-\inf_{\mu\in C}H(\mu|\gamma),-\infty\right\}
	=-\inf_{\mu\in C}H(\mu|\gamma).
	\end{align*}
	Thus, the tail probability in the density comparison is negligible at speed $p_N$ and does not affect the Sanov upper bound.

	Let $O\subseteq\mathcal P(\R)$ be open. If $\inf_{\mu\in O}H(\mu|\gamma)=\infty$, the LDP lower bound is immediate. Otherwise, Sanov's lower bound shows that $\P(\widehat\nu_N\in O)\ge e^{-Mp_N}$ for all sufficiently large $N$ and some $M<\infty$. Since $\P(\|G_N\|_F>R_N)$ is negligible at speed $p_N$, it follows that
	$$
	\frac{\P(\|G_N\|_F>R_N)}{\P(\widehat\nu_N\in O)}\longrightarrow0.
	$$
	The second inequality in Lemma \ref{lem:density_comparison} therefore yields
	\begin{align*}
	\liminf_{N\to\infty}\frac{1}{p_N}\log\P(\nu_N\in O)
	&\ge \liminf_{N\to\infty}\frac{1}{p_N}\log\P(\widehat\nu_N\in O)\\
	&\ge -\inf_{\mu\in O}H(\mu|\gamma).
	\end{align*}
	These bounds prove the claim.
\end{proof}

	\subsection{Proofs of Theorems \ref{thm:MDPforStiefel_main} and \ref{thm:MDPorthogonalgroup}}

	We continue with the proof of Theorem \ref{thm:MDPforStiefel_main}. In parts, we can follow the lines of the proof of \cite[Theorem B]{KabluchkoLarge}. In other parts, the scaling $\beta_N$ causes some technical difficulties that need to be addressed carefully.

	\begin{proof}[Proof of Theorem \ref{thm:MDPforStiefel_main}]
		We start by proving a weak LDP for the block $\beta_N V_{m,N}^{(k)}$, where $V_{m,N}^{(k)}$ denotes the $m\times k$ matrix consisting of the first $k$ columns of $V_{m,N}$. After that, we need to prove exponential tightness to get the full LDP, since unlike in the LDP case, the random variables are no longer contained in a compact set. The claim for the entire matrix then follows in the end by the Dawson-Gärtner Theorem \ref{prop:B:dawson-gaertner}. \\
		
		\noindent Recall that by \cite[Lemma 2.5]{jiang2006many}, the block $V_{m,N}^{(k)}$ has, for $N\ge k+m$, the density $f_N:[-1,1]^{m\times k}\to [0,\infty)$,
		$$
		f_N(A)=\frac{\Gamma_m\left(\frac{N}{2}\right)}{\pi^{\frac{m k}{2}} \Gamma_m\left(\frac{N-k}{2}\right)} \operatorname{det}\left(\operatorname{Id}_{m}-A A^\top\right)^{\frac{N-k-m-1}{2}} \mathbf{1}_{\left\{\left\|A A^\top\right\|_{op}<1\right\}}, \quad A \in[-1,1]^{m \times k}.
		$$
		Therefore, the block $\beta_N V_{m,N}^{(k)}$ has density $\Tilde{f}_N:[-\beta_N,\beta_N]^{m\times k}\to [0,\infty)$ given by
		$$
		\Tilde{f}_N(A) = \beta_N^{-mk}f_N\left(\frac{A}{\beta_N}\right)=\frac{\Gamma_m\left(\frac{N}{2}\right)}{\beta_N^{mk}\pi^{\frac{m k}{2}} \Gamma_m\left(\frac{N-k}{2}\right)} \operatorname{det}\left(\operatorname{Id}_{m}-\frac{1}{\beta_N^2}A A^\top\right)^{\frac{N-k-m-1}{2}} \mathbf{1}_{\left\{\left\|A A^\top\right\|_{op}<\beta_N^2 \right\}},
		$$
		for $A \in[-\beta_N,\beta_N]^{m\times k}$. \\
		
		\noindent\textbf{Step 1: Weak LDP} Define $\mathbb{B}:= \{M\in [-\beta_N,\beta_N]^{m\times k}:\enskip \|MM^\top\|_{op}< \beta_N^2\}$ and its closure by $\overline{\mathbb{B}}:= \{M\in [-\beta_N, \beta_N]^{m\times k}:\enskip \|MM^\top\|_{op}\le \beta_N^2\}$. We want to prove a weak LDP for $\beta_N V_{m,N}^{(k)}$ with speed $\frac{N}{\beta_N^2}$ and rate function 
		$$
		I_k:\R^{m\times k}\to [0,\infty], \quad I_k(V) = \frac{1}{2}\sum_{i=1}^m \|V_i\|^2 = \frac{1}{2}\|V\|_F^2,
		$$ 
		i.e. we need to show the two inequalities according to Proposition \ref{prop:weak-ldp-base}:
		\begin{align}
			& \lim _{r \rightarrow 0} \limsup _{N \rightarrow \infty} \frac{\beta_N^2}{N} \log \mathbb{P}\left(\beta_N V_{m,N}^{(k)} \in B_r(A)\right) \leq -I_k(A),\label{upperboundwaek} \\
			& \lim _{r \rightarrow 0} \liminf _{N \rightarrow \infty} \frac{\beta_N^2}{N} \log \mathbb{P}\left(\beta_N V_{m,N}^{(k)} \in B_r(A)\right) \geq -I_k(A),\label{lowerboundweak}
		\end{align}
		where $B_r(A)$ denotes a Euclidean ball of radius $r>0$ in $\R^{m\times k}$ around $A$ (these sets clearly form a base of the topology). We always have $A\in \overline{\mathbb{B}}$ if $N$ is large enough for any fixed $A\in\R^{m\times k}$. For large $N$,
		\begin{align*}
			&\frac{\beta_N^2}{N} \log \mathbb{P}\left(\beta_N V_{m,N}^{(k)} \in B_r(A)\right) = \frac{\beta_N^2}{N} \log \int_{B_r(A)}\beta_N^{-mk}f_N\left(\frac{D}{\beta_N}\right) \dint D\\
			&= \frac{\beta_N^2}{N} \log \int_{B_r(A)} e^{\frac{N}{\beta_N^2}\frac{\beta_N^2}{N}\log\left(\beta_N^{-mk}f_N\left(\frac{D}{\beta_N}\right)\right)} \dint D \\
	  &\le \frac{\beta_N^2}{N} \log \int_{B_r(A)} e^{\frac{N}{\beta_N^2} \sup_{C\in B_r(A)\cap\mathbb{B}}\frac{\beta_N^2}{N}\log\left(\beta_N^{-mk}f_N\left(\frac{C}{\beta_N}\right)\right)} \dint D\\
			&= \frac{\beta_N^2}{N}\log\left(\vol_{mk}(B_r(A)\cap \mathbb{B})\right) + \sup_{C\in B_r(A)\cap\mathbb{B}}\frac{\beta_N^2}{N}\log\left(\beta_N^{-mk}f_N\left(\frac{C}{\beta_N}\right)\right).
		\end{align*}
		For fixed $A$ and $r$, we have $B_r(A)\subset\mathbb{B}$ for all sufficiently large $N$. Hence, the volume in the first term equals $\vol_{mk}(B_r(A))$. Since $\beta_N^2/N\to0$, this term tends to $0$. Furthermore, we have
		\begin{align*}
			&\frac{\beta_N^2}{N}\log\left(\beta_N^{-mk}f_N\left(\frac{C}{\beta_N}\right)\right)\\
			=& \frac{\beta_N^2}{N} \log\left(\frac{\Gamma_m\left(\frac{N}{2}\right)}{\beta_N^{mk}\pi^{\frac{m k}{2}}\Gamma_m\left(\frac{N-k}{2}\right)} \right) +\frac{\beta_N^2}{N}\frac{N-m-k-1}{2}\log\left(\det\left(\operatorname{Id}_m - \frac{1}{\beta_N^2}CC^\top\right)\right).
		\end{align*}
		Since $m$ and $k$ are fixed, Lemma \ref{lem:gamma_quotient} gives
		\begin{align*}
		&\frac{\beta_N^2}{N}\log\left(\frac{\Gamma_m\left(\frac{N}{2}\right)}{\beta_N^{mk}\pi^{\frac{mk}{2}}\Gamma_m\left(\frac{N-k}{2}\right)}\right)\\
		&=\frac{mk}{2}\frac{\beta_N^2}{N}\log\left(\frac{N}{2\pi\beta_N^2}\right)+O\left(\frac{\beta_N^2}{N^2}\right)\longrightarrow0,
		\end{align*}
		because $x\log(1/x)\to0$ as $x\downarrow0$ and $\beta_N^2/N\to0$.
		For fixed $r>0$, set $M_r:=\sup_{C\in B_r(A)}\|CC^\top\|_{op}<\infty$. For all sufficiently large $N$, we have $B_r(A)\subset\mathbb{B}$. If $\lambda_1(C),\ldots,\lambda_m(C)$ are the eigenvalues of $CC^\top$, then the inequality $|\log(1-x)+x|\le x^2/(2(1-x))$ for $x\in[0,1)$ gives
		\begin{align*}
		&\sup_{C\in B_r(A)}\left|\beta_N^2\log\det\left(\operatorname{Id}_m-\frac{1}{\beta_N^2}CC^\top\right)+\|C\|_F^2\right|\\
		&\le \frac{1}{2\beta_N^2(1-M_r/\beta_N^2)}\sup_{C\in B_r(A)}\operatorname{Tr}\left((CC^\top)^2\right)\longrightarrow0.
		\end{align*}
		Since $(N-m-k-1)/N\to1$, it follows that
		\begin{align*}
		&\limsup_{N\to\infty}\sup_{C\in B_r(A)}\frac{\beta_N^2}{N}\frac{N-m-k-1}{2}\log\det\left(\operatorname{Id}_m-\frac{1}{\beta_N^2}CC^\top\right)\\
		&=-\frac{1}{2}\inf_{C\in B_r(A)}\|C\|_F^2.
		\end{align*}
		Consequently,
		\begin{align*}
		\lim _{r \rightarrow 0} \limsup _{N \rightarrow \infty} \frac{\beta_N^2}{N} \log \mathbb{P}\left(\beta_N V_{m,N}^{(k)} \in B_r(A)\right) \le - \frac{1}{2}\lim _{r \rightarrow 0}\inf_{C\in B_r(A)}\|C\|_F^2 = -\frac{1}{2}\|A\|_F^2,
		\end{align*}
		by continuity in $C$. The upper bound \eqref{upperboundwaek} is satisfied. \\

		For the lower bound, the density integral satisfies, for all sufficiently large $N$,
		\begin{align*}
		&\frac{\beta_N^2}{N}\log \mathbb{P}\left(\beta_N V_{m,N}^{(k)} \in B_r(A)\right)\\
		&\ge \frac{\beta_N^2}{N}\log\vol_{mk}(B_r(A))
		+\inf_{C\in B_r(A)}\frac{\beta_N^2}{N}\log\left(\beta_N^{-mk}f_N\left(\frac{C}{\beta_N}\right)\right).
		\end{align*}
		The volume term and the normalizing constant tend to $0$. The uniform determinant estimate therefore gives
		\begin{align*}
		\liminf_{N\to\infty}\frac{\beta_N^2}{N}\log \mathbb{P}\left(\beta_N V_{m,N}^{(k)} \in B_r(A)\right)
		&\ge -\frac{1}{2}\sup_{C\in B_r(A)}\|C\|_F^2.
		\end{align*}
		Letting $r\to0$ and using continuity proves the lower bound \eqref{lowerboundweak}.
		Finally, $I_k$ is continuous and therefore lower semicontinuous. Also 
		$$
		\{A\in\R^{m\times k}:\enskip I_k(A)\le y\}
		$$ 
		is the closed Frobenius ball of radius $\sqrt{2y}$ and is therefore compact, which makes $I_k$ a good rate function.\\
		
		\noindent We proved that the sequence $\beta_N V_{m,N}^{(k)}$, for uniform random $V_{m,N}$ in $\mathbb{V}_{m,N}$, satisfies a weak LDP with speed $\frac{N}{\beta_N^2}$ and good rate function $I_k:\R^{m\times k}\to [0,\infty]$, $I_k(A) = \frac{1}{2}\sum_{i=1}^m \|A_i\|^2$. \\
		
		\noindent\textbf{Step 2: Exponential tightness} 
		We need to show that for any $\epsilon>0$, there exists a compact set $K\subset \R^{m\times k}$ such that
		$$
		\limsup_{N\to\infty}\frac{\beta_N^2}{N}\log(\P(\beta_N V_{m,N}^{(k)}\in K^c)) \le -\epsilon.
		$$
		Represent $V_{m,N}$ by the first $m$ rows of a Haar-distributed matrix in $\O(N)$. For fixed $i$ and $j$, the map $U\mapsto U_{i,j}$ is $1$-Lipschitz with respect to the Hilbert-Schmidt distance. For all sufficiently large $N$, multiplication on the left by a diagonal matrix that changes the signs of rows $i$ and $r$, where $r\ne i$, preserves each component of $\O(N)$ and changes the sign of $U_{i,j}$. Thus, the coordinate has mean zero on each component. Applying Proposition 10 in \cite{meckes2017rates} to the two components gives constants $c,C>0$ such that
		$$
		\P(|(V_{m,N})_{i,j}|>t)\le C e^{-cNt^2}
		$$
		for all $t>0$ and all sufficiently large $N$. For $L>0$, set $K=[-L,L]^{m\times k}$. A union bound gives
		$$
		\P(\beta_N V_{m,N}^{(k)}\in K^c)
		\le Cmk\exp\left(-\frac{cNL^2}{\beta_N^2}\right).
		$$
		Since $\beta_N^2/N\to0$, it follows that
		$$
		\limsup_{N\to\infty}\frac{\beta_N^2}{N}\log\P(\beta_N V_{m,N}^{(k)}\in K^c)\le -cL^2.
		$$
		Choosing $L$ such that $cL^2\ge\epsilon$ proves exponential tightness.
		This concludes the proof of the exponential tightness. We proved that the sequence $\beta_N V_{m,N}^{(k)}$ for uniform random $V_{m,N}$ in $\mathbb{V}_{m,N}$, satisfies an LDP on $\R^{m\times k}$ with speed $\frac{N}{\beta_N^2}$ and good rate function $I_k:\R^{m\times k}\to [0,\infty]$, $I_k(A) = \frac{1}{2}\sum_{i=1}^m \|A_i\|^2$. \\
		
		\noindent\textbf{Step 3: Projective limit} We want to uplift the LDP to the entire sequence $\beta_N V_{m,N}$. The projective system can be defined by $\left({Y}_j, p_{i j}\right)_{i \leq j}$, with ${Y}_k =\R^{m\times k}$, and define for $i\le j$, $p_{ij}(V) = V^{(i)}$ (i.e. it removes the columns $i+1,...,j$). This clearly forms a projective system. Then, the projective limit is ${Y} := \R^{m\times\infty}$, with projection maps $p_k(V) = V^{(k)}$. We proved that the sequence $p_k(\beta_N V_{m,N})= \beta_N V_{m,N}^{(k)}$ satisfies an LDP with speed $\frac{N}{\beta_N^2}$ and good rate function $I_k$ for every fixed $k$. By the Dawson-Gärtner Theorem \ref{prop:B:dawson-gaertner}, the sequence $\beta_N V_{m,N}$ satisfies an LDP with speed $\frac{N}{\beta_N^2}$ and good rate function $I:\R^{m\times \infty}\to [0,\infty]$, $I(V) = \sup_{k}I_k(p_k(V))$. Since $I_k(p_k(V))=\frac{1}{2}\sum_{i=1}^m\sum_{j=1}^k v_{i,j}^2$, taking the supremum over $k$ gives $I(V) = \frac{1}{2}\sum_{i=1}^m \|V_i\|^2$. The convexity of this rate function is clear. This concludes the proof. 
		\end{proof}

		\begin{proof}[Proof of Theorem \ref{thm:MDPorthogonalgroup}]	
		Recall that $\mathbb{V}_{N,N} = \O(N)$. We define the projective system $\left({Y}_j, p_{i j}\right)_{i \leq j}$, by ${Y}_m =\R^{m\times \infty}$, and define for $i\le j$, $p_{ij}(V) = V^{(i, \infty)}$ (i.e. it removes the rows $i+1,...,j$). The projective limit is ${Y} = \R^{\infty\times\infty}$, with projection maps $p_k(V) = V^{(k,\infty)}$. Let $A_N$ be chosen uniformly from $\O(N)$ for each $N$. Then $A_N^{(k,\infty)}$ is distributed uniformly in $\mathbb{V}_{k,N}$. By Theorem \ref{thm:MDPforStiefel_main}, the sequence $p_m(\beta_N A_N) = \beta_N A_N^{(m,\infty)}$ satisfies an LDP with speed $\frac{N}{\beta_N^2}$ and good convex rate function $I_m:\R^{m\times\infty}\to[0,\infty]$, $I_m(A) = \frac{1}{2}\sum_{i=1}^m \|A_i\|^2$. By the Dawson-Gärtner theorem \ref{prop:B:dawson-gaertner}, the sequence $\beta_N A_N$ satisfies an LDP with speed $\frac{N}{\beta_N^2}$ and good convex rate function $J(A) = \sup_{m}I_m(p_m(A)) = \frac{1}{2}\sum_{i=1}^\infty \|A_i\|^2$, which proves the claim. 
		\end{proof}

\section*{Acknowledgments}

The author is supported by the DFG project \emph{Limit theorems for the volume of random projections of $\ell_p$-balls} (project number 516672205). The author also thanks Joscha Prochno and Christoph Th\"ale for a careful reading of a preliminary version of this paper and for many helpful suggestions.

\end{document}